\documentclass[12pt]{article}
\usepackage{amssymb}

\def\aa{{\mathcal A}}
\def\al{\alpha}
\def\be{\beta}
\def\ff{{\mathcal F}}
\def\forces{{\;\Vdash\;}}
\def\ga{\gamma}
\def\goth{\mathfrak}
\def\om{\omega}
\def\poset{{\mathbb P}}
\def\prodaa{\prod_{n<\om}\aa_n}
\def\proof{\par\noindent proof:\par} 
\def\qed{\nopagebreak\par\noindent\nopagebreak$\Box$\par}
\def\reals{{\mathbb R}}
\def\res{\mathord{\upharpoonright}}
\def\rmand{\mbox{ and }}

\def\sm{\setminus}
\def\su{\subseteq}

\newtheorem{theorem}{Theorem} 
\newtheorem{prop}[theorem]{Proposition}

\begin{document}

\begin{center}
The onto mapping property of Sierpinski
\end{center}

\begin{flushright}
A. Miller \\
July 2014
\end{flushright}

\bigskip

\noindent Define

(*) There exists $(\phi_n:\om_1\to \om_1:n<\om)$ such that
for every $I\in [\om_1]^{\om_1}$ there exists $n$ such
that $\phi_n(I)=\om_1$.

\bigskip

This is roughly what Sierpinski \cite{sier} refers to as $P_3$ but I think
he brings $\reals$ into it.  I don't know French so I cannot say for sure
what he says but I think he proves that (*) follows from
the continuum hypothesis.  
Here we show that the existence of a Luzin set implies (*) and 
(*) implies that there exists a nonmeager set of reals of size $\om_1$.
We also show that it is relatively consistent that (*) holds
but there is no Luzin set.  All the other properties in this paper,
(**), (S*), (S**), (B*) are shown to be equivalent to (*).

\bigskip

\begin{prop}(Sierpinski \cite{sier})
CH implies (*).
\end{prop}

\proof

Let $\om_1^\om=\bigcup_{\al<\om_1}\ff_\al$ 
where the $\ff_\al$ are countable and increasing.
For each $\al$ construct $(\phi_n(\al):n<\om)$ so that for
every $g\in \ff_\al$ there is a some $n$ such that $\phi_n(\al)=g(n)$.

Now suppose $I\su\om_1$.  If no $\phi_n$ maps $I$ onto $\om_1$, 
then there
exists $g\in\om_1^\om$ such that $g(n)\notin\phi_n(I)$ for every $n$.
If $g\in \ff_{\al_0}$, then $\al\notin I$ for every $\al\geq\al_0$.
This is because $g\in \ff_\al$ and so for some $n$ $g(n)=\phi_n(\al)$
and since $g(n)\notin \phi_n(I)$ we have $\al\notin I$.
\qed

\bigskip \noindent Define

(**) There exists $(g_\al:\om\to\om_1:\al<\om_1)$ such that for every
$g:\om\to\om_1$ for all but countably many $\al$ there are infinitely
many $n$ with $g(n)=g_\al(n)$.

\begin{prop}
(**) iff (*).
\end{prop}

\proof

To see (**) implies (*) let $\phi_n(\al)=g_\al(n)$. Then the proof
of the first proposition goes thru. 

On the other hand suppose
$(\phi_n:\om_1\to \om_1:n<\om)$ witnesses (*).  First note
that for any 
$I\in [\om_1]^{\om_1}$ there are infinitely many $n$ such
that $\phi_n(I)=\om_1$.  This is because if there are only
finitely many $n$ we could cut down $I$ in finitely many steps
so that there were no
$n$ with $\phi_n(I)=\om_1$.

Now define $g_\al\in\om_1^\om$ by 
$g_\al(n)=\phi_n(\al)$.  These witness (**).
Given any $g:\om\to\om_1$ if there is
an uncountable $I\su\om_1$ and $N<\om$ such that for every $\al\in I$
we have $g(n)\neq g_\al(n)$ for all $n>N$ then this means
that $g(n)\notin \phi_n(I)$ and for all $n>N$ and so (*) fails.

\qed

Obviously (**) is false if ${\goth b}>\om_1$
so (*) is not provable just from ZFC.

\begin{prop}
It is relatively consistent with any cardinal arithmetic that 
(*) is true and ${\goth b}={\goth d}=\om_1$.
\end{prop}

\proof

Start with any $M$ a countable transitive model of ZFC.
Our final model is $M[g_\al,f_\al:\al<\om_1]$
where each $g_\al:\om\to\al$ is generic with respect to the poset of
finite partial functions from $\om$ to $\al$ and $f_\be\in\om^\om$ is
Hechler real over $M[g_\al,f_\al:\al<\be]$.  The $\om_1$-sequence
is obtained by finite support ccc forcing.
By ccc for any $g\in \om_1^\om\cap M[g_\al,f_\al:\al<\om_1]$
there will be $\al_0<\om_1$ such 
that $\al_0$ bounds the range of $g$ and
$g\in M[g_\al,f_\al:\al<\al_0]$.  It follows by product genericity that
for every $\al\geq\al_0$ there are infinitely many $n$ such that
$g(n)=g_\al(n)$.  The Hechler sequence $f_\al$ for $\al<\om_1$ shows
that ${\goth d}=\om_1$.

\qed

With a little more work we will prove that (*) follows
from the existence of a Luzin set (Prop \ref{luzin}).
We will also show that (*) implies there is a nonmeager set
of reals of size $\om_1$ (Prop \ref{nonmeag}) and so in the
random real model (*) fails and ${\goth b}={\goth d}=\om_1$.

\bigskip

\bigskip

\bigskip\noindent
Actually I think Sierpinski 
considers what appears to be a stronger version:

\bigskip\noindent Define

(S*)  There exists $(\phi_n:\om_1\to \om_1:n<\om)$ such that
for every $I\in [\om_1]^{\om_1}$ for all but finitely many $n$
$\;\;\;\phi_n(I)=\om_1$.

\bigskip
Surprisingly (S*) is equivalent to (*).

\begin{prop}\label{sstar}
(S*) iff (*).
\end{prop}

\proof

We show (**) implies (S*).

Let $a_0=1$ and $a_{n+1}=1+\sum_{i\leq n} a_i$.
Let $$\aa_n=\{u \;\;|\;\; \exists D\in [\om_1]^{a_n}\;\; u:D\to\om_1\}
\rmand \prodaa=\{g  \;\;|\;\;  \forall n \; g(n)\in \aa_n\}$$
Since each $\aa_n$ has cardinality $\om_1$ from (**) we
get $(g_\al\in \prodaa:\al<\om_1)$ such that for every 
$g\in\prodaa$ for all but countably many $\al$ there are
infinitely many $n$ such that $g(n)=g_\al(n)$.  For each 
$\al<\om_1$ define $h_\al:\om\to\om_1$ so that if $g_\al(n)=u_n:A_n\to\om_1$
for every $n$ then 
$$h_\al\res (A_n\sm \cup_{i<n}A_i) = u_n\res (A_n\sm \cup_{i<n}A_i)$$
Since $|A_k|=a_k$ the sets $A_n\sm \cup_{i<n}A_i$ are nonempty.
We claim that
the $h_\al$ have the following property:

\bigskip\noindent Define \par
(S**) For any $X\in [\om]^\om$ and $h:X\to\om_1$ for all but countably many
$\al$ there are infinitely many $n\in X$ with $h(n)=h_\al(n)$.

\bigskip
It is enough to see there is at least one $n\in X$ with $h(n)=h_\al(n)$.
Otherwise if there were only finitely many $n$ for uncountably many $\al$
we could throw out from $X$ a fixed finite set for uncountably many $\al$ and
get a contradiction.

Let $X=\{x_n: n<\om\}$ listing $X$ in increasing order.
Define $g\in\prodaa$ by $g(n)= h\res\{x_i:i<a_n\}$.  Now suppose
$g_\al(n)=g(n)$.  This means that if $g_\al(n)=u_n:A_n\to\om_1$, then
$A_n=\{x_i:i<a_n\}$ and $u_n = h\res A_n$.  But since
$A_n\sm \cup_{i<n}A_i$ is nonempty we get that $h_\al(x)=h(x)$ for
some $x\in X$.

Now define $\phi_n(\al)=h_\al(n)$.  This has the required property (S*).
Given $I$ uncountable let $X$ be the $n\in \om$ with $\phi_n(I)\neq \om_1$.
If $X$ is infinite we would get $h:X\to\om_1$ such that
$h(n)\notin \phi_n(I)$ for all $n\in X$. 
But this means that for all $\al\in I$
and $n\in X$ that $h(n)\neq h_\al(n)$ which contradicts (S**).

\qed

This is related to results in Bartoszynski \cite{barto}.

\bigskip
Bagemihl-Sprinkle \cite{bag} say that Sierpinski states CH implies (S*) but
only proves (*).  They give a proof from CH of a
seemingly stronger version:

\bigskip\noindent Define \par
(B*)  There exists $(\phi_n:\om_1\to \om_1:n<\om)$ such that
for every $I\in [\om_1]^{\om_1}$ for all but finitely many $n$
for all $\be<\om_1$ there are uncountably many
$\al\in I$ with $\phi_n(\al)=\be$, i.e., not only is
$\phi_n(I)=\om_1$ but it is uncountable-to-one.

\begin{prop}
(S*) iff (B*)
\end{prop}

\proof
Let $\pi:\om_1\to\om_1$ be uncountable to one, i.e.,
for all $\be<\om_1$ there are uncountably many
$\al<\om_1$ with $\pi(\al)=\be$.  If $(\phi_n:\om_1\to \om_1:n<\om)$ witness
(S*) then $(\pi\circ\phi_n:\om_1\to \om_1:n<\om)$ satisfies (B*).
\qed

\begin{prop}\label{luzin}
If there is a Luzin set, then (*) is true.
\end{prop}
\proof

We prove (**).
Suppose $\{g_\al:\om\to\om\;:\;\al<\om_1\}$ is a Luzin set,
then it satisfies that
for every $k:\om\to\om$ for all but countably many $\al<\om_1$
there are infinitely many $n$ such that $k(n)=g_\al(n)$.

There is a sequence
$(f_\al:\al\to\om:\om\leq \al<\om_1)$ of one-to-one
functions which is coherent: for $\al<\be\;\; f_\be\res\al=^* f_\al$,
i.e., $f_\be(\ga)=f_\al(\ga)$ 
for all but finitely many $\ga<\al$.  This is the
construction of an Aronszajn tree which appears in
the first edition of Kunen's set theory book \cite{kunen}.

Let $\hat{g}_\al:\om\to\al$ be any map which extends $f_\al^{-1}\circ g_\al$.
We claim that for any $k:\om\to\om_1$ which is one-to-one that
for all but countably many $\al$ there are infinitely many $n$ with
$\hat{g}_\al(n)=k(n)$.   To see this suppose $k:\om\to\be$ is one-to-one
and let $\hat{k}=f_\be\circ k$ which maps $\om$ to $\om$. Then for
some $\al_0>\be$ for all $\al\ge\al_0$
there will be infinitely
many $n$ with $g_\al(n)=\hat{k}(n)$.  This means
that $g_\al(n)=f_\be(k(n))$.  Since $k$ is one-to-one, there
will be infinitely many such $n$ where $f_\be(k(n))=f_\al(k(n))$.
But $g_\al(n)=f_\al(k(n))$ implies $\hat{g}_\al(n)=k(n)$.

To get rid of the requirement that $k$ be one-to-one,
let $j:\om_1\times\om \to \om_1$ be a bijection and
$\pi:\om_1\to \om_1$ be projection onto first coordinate, i.e.,
$\pi(j(\al,n))=\al$.   Define $h_\al(n)=\pi(\hat{g}_\al(n))$. 
Given any $k:\om\to\om_1$ define $\hat{k}(n)=j(k(n),n)$.
Then since $\hat{k}$ is one-to-one for all but countably
many $\al$ there will be infinitely many $n$ with
$\hat{g}_\al(n)=\hat{k}(n)$.  But this implies 
$$h_\al(n)=\pi(\hat{g}_\al(n))=\pi(\hat{k}(n))=k(n)$$
Hence $(h_\al:\al<\om_1)$ satisfies (**).

\qed

\begin{prop}\label{nonmeag}
Suppose (*), then there exists $(x_{\al,\be}\in 2^\om:\al,\be<\om_1)$
such that for every dense open $D\su 2^\om$  there exists
$\al_0<\om_1$ such that for every $\al\geq\al_0$ there is
a $\be_\al<\om_1$ such that $x_{\al,\be}\in D$ for every $\be\geq\be_\al$.
\end{prop}
\proof
We use that there are
$\{h_\al:\om\to\om\;:\;\al<\om_1\}$ with the property that
for every $X\in[\om]^\om$ and $h:\om\to\om$ for all but countably
many $\al$ there are infinitely many $n\in X$ with $h(n)=h_\al(n)$
(see (S**) in the proof of Prop \ref{sstar}).
This implies that there exists 
$(X_\al\in [\om]^\om:\al<\om_1)$ such that for every $Y\in [\om]^\om$
for all but countably many $\al$ there are infinitely many
$x\in X_\al$ such that $|Y\cap [x,x^+)|\geq 2$ where $x^+$ is
the least element of $X_\al$ greater than $x$.  Fix $\al$ and
enumerate $X_\al=\{k_n:n<\om\}$ in strict increasing order.
Define
$$P_\al=\{g:\om\to FIN(\om,2)\;:\;\forall n\;\; g(n)\in 2^{[k_n,k_{n+1})}\}$$
By (S**) there exists $g_{\al,\be}\in P_\al$ for $\be<\om_1$ with the
property that for any $h$ in $P_\al$ and infinite $Y\su\om$ for
all but countably many $\be$ there are infinitely many $n\in Y$ with
$h(n)=g_{\al,\be}(n)$.   Define $x_{\al,\be}\in 2^\om$ by
$x_{\al,\be}(m)=g_{\al,\be}(n)(m)$ where $n$ is the unique integer
with $k_n\leq m<k_{n+1}$.  Equivalently $x_{\al,\be}=\bigcup_ng_{\al,\be}(n)$.
(Without loss we may assume $k_0=0\in X_\al$.)

Given $D\su 2^\om$ dense open let $\hat{D}\su 2^{<\om}$ be the
set of all $s$ with $[s]\su D$.  Construct an infinite $Z\su\om$ so
that for every $z\in Z$ there exists $t\in 2^{<\om}$ with
$|t|\leq z^+-z$ such that for every $s\in 2^{<\om}$ with
$|s|\leq z$ we have $s\, t\in \hat{D}$ where $s\,t$ is
the concatenation of $s$ with $t$.  By construction
there exists $\al_0$ so that for every $\al\geq\al_0$ the
there are infinitely many $x\in X_\al$ with $|[x,x^+)\cap Z|\geq 2$.

Fix $\al\geq \al_0$ and as above $X_\al=\{k_n:n<\om\}$.  Let
 $$Y=\{n\;:\; |[k_n,k_{n+1})\cap Z|\geq 2.$$
Note that by the definition
of $Y$ there is a $h\in P_\al$ with the property that for
every $n\in Y$ for every $s\in 2^{k_n}$ we have $s\cup h(n)\in \hat{D}$.
For some $\be_\al$ for every $\be\geq\be_\al$ there are infinitely
many $n\in Y$ with $h(n)=g_{\al,\be}(n)$ and so $x_{\al,\be}\in D$.
\qed

This is similar to the argument of Miller \cite{char}.
Obviously the set of $x_{\al,\be}$ in Prop \ref{nonmeag} is nonmeager.
Although it seems a little bit like a Luzin set, it isn't.

\begin{prop} \label{superperf}
In the superperfect tree model (*) holds but there is no Luzin set.
\end{prop}

\proof This is the countable support iteration of length $\om_2$ of 
superperfect tree forcing\footnote{ So called Miller forcing. I also called it
rational perfect set forcing. } over a ground model of CH. The fact that there
is no Luzin set in this model is due to  Judah and Shelah \cite{judah}.  They
also show that the set of ground model reals is not meager. We first do the
argument for a single superperfect real even though it is not needed but
it is easy and allows us to show the rest of the argument.  Then we
quote known results to cover the countable support iteration of length $\om_2$.

For $T$ a subtree of $\om^{<\om}$, a node $s\in T$ is a splitting
node iff $sn\in T$ for infinitely many $n<\om$.
A tree $T\su \om^{<\om}$ is superperfect iff the splitting nodes
of $T$ are dense in the tree $T$. The poset $\poset$ is the partial
order of superperfect trees.

\bigskip\noindent {\bf One Step Lemma}.
Suppose $p\in\poset$, $\al<\om_1$, $\tau$ is a $\poset$-name such that
$p\forces \tau:\om\to\al$, and $X\in [\om]^\om$.  Then
there exists $f:X\to\al$ and $q\leq p$ such that
$$q\forces \exists^\infty n\in \check{X} \; \check{f}(n)=\tau(n)$$

\proof
To prove this lemma,
let $\{x_s \;:\;s\in p\}$ be a one-to-one enumeration of $X$.
By standard fusion arguments construct $q\leq p$ and $f$ such that
for every split node $s\in q$ and $sn\in q$ we have that
$$q_{sn}\forces \check{f}(x_{s_n})=\tau(x_{s_n})$$
\qed

\bigskip

Now we show that we can construct a witness to (**) which remains
one after forcing once with $\poset$.  Let $X_\al\in [\om]^\om$
for $\al<\om_1$ be pairwise disjoint.  Let $\{(p_\al,\tau_\al):\al<\om_1\}$
list all pairs of $(p,\tau)$ such that $p\in\poset$ and
$\tau$ is a canonical name such that $p\forces \tau:\om\to\om_1$.
Apply the One Step Lemma to get $q_\al\leq p_\al$ and
$f_\al:X_\al\to\om_1$ such that 

$$q_\al\forces \exists^\infty n\in \check{X_\al} \;\; 
\check{f_\al}(n)=\tau_\al(n)$$

Now construct $g_\al:\om\to\om_1$ such that for every $\be<\al$
$\;\;g_\al\res X_\be=^*f_\be$.  (To see how to do this
let $\{\be_n:n<\om\}$ be a one-to-one enumeration of $\al$.
Put $Z_n=X_{\be_n}\sm \bigcup_{k<n}X_{\be_k}$ and
$g_\al=\bigcup_{n<\om}f_{\be_n}\res Z_n$.)

We claim after forcing with $\poset$ that
$(g_\al:\al<\om_1)$ satisfies $(**)$.   Suppose $p\forces \tau:\om\to\om_1$.
We may find $p_\al\leq p$ and $\tau_\al$ such that
$p_\al\forces \tau=\tau_\al$.  By construction 
$$q_\al \forces \exists^\infty n\in X_\al\;\; f_\al(n)=\tau_\al(n)$$
Since for any $\ga>\al$ we have that $g_\ga\res X_\al=^*f_\al$ 
we are done.

\bigskip
The next step is to generalize the One Step Lemma to $\poset_{\om_2}$
by using a result of Judah and Shelah \cite{judah}. They showed that after
forcing with $\poset_{\om_2}$
the set of ground model reals,  $M\cap \om^\om$, is nonmeager. Hence
for any $X\in[\om]^\om\cap M$ and $\al<\om_1$ we
have that $M\cap \al^X$ is nonmeager. 
Thus for any $k:X\to\om$ in the generic extension $M[G]$ there must be
$f:X\to\om$ in $M$ such that $f(n)=k(n)$ for infinitely many $n\in X$.
This is because the set 
$$\{f\in \al^X:\forall^\infty n\in X\;\; f(n)\neq k(n)\}$$
is meager.
Hence the Lemma holds for $\poset_{\om_2}$, i.e.,
for any $\tau$, $X\in[\om]^\om$, $\al<\om_1$ and $p\in \poset_{\om_2}$
such that $p\forces \tau:\om\to\al$
there is $f\in \al^X$ and $q\leq p$ such that 
$$q\forces \exists^\infty n\in \check{X} \; \check{f}(n)=\tau(n).$$

Superperfect tree forcing is Souslin; Goldstern and Judah \cite{gold} give
the argument in detail for Laver forcing.
An earlier paper of Judah and Shelah \cite{sous} shows that
every real in the $\om_2$ length iteration of Souslin posets
is added by a sub-iteration of countable length.
Hence for any $G_{\om_2}$ which is $\poset_{\om_2}$ generic over $M$ and
$k\in 2^\om\cap M[G_{\om_2}]$
there exists $\al<\om_1$ and $H_\al\in M[G_{\om_2}]$ which is $\poset_\al$-generic
over $M$ with $k\in M[H_\al]$.
Judah and Shelah \cite{sous} do this in detail for the iteration
of Mathias forcing but
it would also be true for the iteration of superperfect tree forcing.
Hence we only need worry about pairs of conditions and names for 
$\poset_{\al}$ for $\al<\om_1$.  Up to forcing equivalence there are
only $\om_1$ of them.

This proves Proposition \ref{superperf}.
\qed

\bigskip\noindent Does the existence of a nonmeager set of reals of size
$\om_1$ imply (*)?

\bigskip

This paper was motivated by a result in an earlier
version of A.Medini \cite{med} which
showed that (*) implies that there is an uncountable $X\su 2^\om$ with 
the Grinzing property: for every uncountable
$Y\su X$ there is an uncountable family of uncountable
subsets of $Y$ with pairwise disjoint closures in $2^\om$.
To do this Medini used a result from Miller \cite{onto}.
This has been superceded by a proof in ZFC of an uncountable $X\su 2^\om$
with the Grinzing property.

\bigskip

\begin{flushleft}
 Arnold W. Miller \\
 miller@math.wisc.edu \\
 http://www.math.wisc.edu/$\sim$miller\\
 University of Wisconsin-Madison \\
 Department of Mathematics, Van Vleck Hall \\
 480 Lincoln Drive \\
 Madison, Wisconsin 53706-1388 \\
\end{flushleft}

\end{document}